\magnification=1000

\def\cqfd{\unskip\kern 6pt\penalty 500
\raise 0pt\hbox{\vrule\vbox to5pt{\hrule width 5pt
\vfill\hrule}\vrule}\par}

{\parindent=0mm{\footnote{}{2000 {Mathematics Subject Classification}.
Primary 37E30 \obeylines
}}} 
 \centerline{\bf FREE LINES FOR HOMEOMORPHISMS OF THE OPEN ANNULUS}
\bigskip
\centerline{LUCIEN GUILLOU }
\bigskip
\bigskip
\noindent {\bf Abstract}. Let  $ H $ be a homeomorphism of the open annulus $S^1
\times {\bf R}$ isotopic to the identity and let  $h$ be a lift of $H$ to
the universal cover ${\bf R} \times {\bf R}$ without fixed point. Then we show that $h$ admits a Brouwer line which is a lift of a properly imbedded line joining one end to the other in the annulus or $H$ admits a free essential simple closed curve.
  \bigskip
 \bigskip

\noindent {\bf 0. Introduction.}

\bigskip

Our main result in this note is the following:

\medskip

{\bf 0.1 Theorem}. {\it Let $H$ be a homeomorphism of the open annulus $S^1 \times {\bf R}$ isotopic to the identity which admits a lift $h$ to ${\bf R}^2$ without fixed point. Then

there exists an essential simple closed curve in the annulus free under $H$ (and therefore it lifts to a Brouwer line for $h$)

or there exists a properly imbedded line in the annulus joining both ends which lifts to a Brouwer line for $h$.}

\medskip

A natural question, in case $H$ is fixed point free, is to ask, in the second case, for a  properly imbedded line in the annulus joining both ends and {\bf free} for $H$. It has been given a positive answer in many cases recently by  P. B\'eguin, S. Crovisier and F. Le Roux [BCL].

\bigskip

{\bf Note}. A slightly weaker form of this theorem was communicated to a few people since 1994 but not published partly because I was unable to answer the preceding question. In his 2000 thesis, Sauzet [S] reproved in his own way the exact statement above but this thesis stayed also unpublished. Recently, Addas-Zanatta [A] and B\'eguin-Crovisier-Le Roux [BCL] used this result as a step towards interesting dynamical results, so we decided to publish these notes. Our original proof relied on hard to read original work of K\'er\'ekjart\'o [K2], we tried here to give a proof starting from first principles.

Note further that, as remarked in [BCL], this theorem, for area preserving homeomorphisms, follows quickly from the much harder recent ``version feuillet\'e \'equivariante du th\'eor\`eme de Brouwer'' of Le Calvez [L].

\bigskip

For comparison recall that in [G] we proved that {\it if $H$ is a fixed point free orientation preserving homeomorphism of the closed annulus, then there exists an essential simple closed curve in the annulus free under $H$ or there exists an essential simple arc in the annulus free under $H$} generalizing the celebrated Poincar\'e-Birkhoff theorem. Part of the arguments given here can easily be adapted to prove this result too. Instead, these arguments will be used to prove:

\bigskip

{\bf 0.2 Theorem}. {\it Let $H$ be a homeomorphism of the torus ${\bf T}^2$ preserving orientation which admits a lift to ${\bf R}^2$ without fixed point. Then there exists an essential simple closed curve on ${\bf T}^2$ which lifts to a Brouwer line of $h$}.

Note that in the case of the torus, we cannot ask for an essential simple curve in ${\bf T}^2$ free under $H$: Bestvina and Handel have given an example of a fixed point free, homotopic to the identity, homeomorphism of the torus without any free essential simple closed curve [BH].

\medskip

Finally, in section 5, we show how Theorem 0.1 can be used to solve a conjecture of Winkelnkemper about area preserving homeomorphisms of the 2-sphere with exactly two stable fixed points which generalises the classical Poincar\'e-Birkhoff Theorem on the closed annulus.

\bigskip

{\bf Acknowledgements}. I would like to thank F. Le Roux for useful comments on a preliminary version of this paper.

\vskip 2cm

\noindent{\bf 1. Basic facts about Brouwer homeomorphisms.}

\bigskip

In this paragraph $h$ will be a fixed point free orientation preserving homeomorphism of the
plane ${\bf R}^2$. Such an homeomorphism is called a {\bf Brouwer homeomorphism}. 

The most basic fact about these homeomorphisms is that they automatically no periodic points or even recurrent points. More generally we have the following lemma of J. Franks [F] for which we recall the following definition.

\bigskip

{\bf 1.1 Definition}. Let $f$ be a homeomorphism of ${\bf R}^2$. A subset $E$ of ${\bf R}^2$ is {\bf free} (for $f$) if $E\bigcap f(E)=\emptyset$. It is {\bf critical} if int$E$ is free but $f(E) \bigcap E \not = \emptyset$.

\bigskip

{\bf 1.2 Lemma}. {\it Let $(B_i)_{i=1, \ldots , k}$, $(k \geq 1)$ be a sequence of open discs in ${\bf R}^2$ and $f$ an orientation preserving homeomorphism of ${\bf R}^2$ such that:

 1. $B_i \bigcap B_j = \emptyset$ for all $i \not = j, 1 \leq i,j \leq k$;

2. each $B_i$ is free for $f$ for each $1 \leq i \leq k$;

3. there exist integers $n_i >0$ such that $f^{n_i}(B_i) \bigcap B_{i+1} \not = \emptyset$ for each $i \in {\bf N}$ where $i$ is taken mod\ $k$.

Then there exists a simple closed curve of index $1$ for $f$ and, therefore, a fixed point for $f$.}

\bigskip

We fix now a Brouwer homeomorphism $h$. As a consequence of Franks' Lemma, every point in the plane is wandering under $h$.

\bigskip

{\bf 1.3 Definition}. A simple closed arc $\alpha \subset {\bf R}^2$ with endpoints $p$ and
$h(p)$ for some $p\in {\bf R}^2$ is a {\bf translation arc} if $\alpha \bigcap h(\alpha )=\{ h(p)\}$. A translation arc will always be supposed to be oriented from $p$ to $h(p)$.

\bigskip

{\bf 1.4 Proposition}. {\it Every point of ${\bf R}^2$ is contained in a translation arc. If $\alpha $ is a translation arc, the {\bf trajectory} $\bigcup_{n\in {\bf Z}}h^n (\alpha )$  is a submanifold of ${\bf R}^2$ homeomorphic to ${\bf R}$.}

\bigskip

Such a trajectory $l$ is not necessarily closed as a subset of ${\bf R}^2$ but, being a submanifold, $\bar l \setminus l$ is closed and we have:

\bigskip

{\bf 1.5 Proposition}. {\it A trajectory $l$ divides the connected component of ${\bf R}^2 \setminus (\bar l \setminus l)$ containing $l$ into exactly two components.}

\bigskip

We refer to [G] for a proof of the two preceding propositions.

\vskip 2cm

\noindent{\bf 2. Critical discs and Brouwer plane translation theorem}.

\bigskip

The content of this section may be seen as a variation on [T].

\medskip

Let still $h$ be a Brouwer homeomorphism.

\bigskip

A disc will always be a {\bf closed} disc in that paper, otherwise we say ``open disc''.

\bigskip

{\bf 2.1 Lemma}. {\it If the disk $D$ is free or critical, then $D \bigcap h^n(D) = \emptyset$ for every $n \in {\bf Z} \setminus \{0, \pm 1\}$.}

\medskip

{\it Proof}. Let us suppose by contradiction that there exists a point $x \in D$ such that $h^n(x) \in D$ for some $n>1$. Any arc from $x$ to $h^n(x)$ contained in int$D$ except, eventually, for its endpoints, would be free since $h$ has no periodic point. Then, one could find a open disc, neighborhhod of this arc which contradicts Franks' Lemma (with $i=1$).

\bigskip

{\bf 2.2 Lemma}. {\it Let $D$ be a critical disc. Then the boundary of $D$ decomposes into two disjoint arcs  (eventually degenerated, that is reduced to a point) containing respectively $h^{-1}(D) \bigcap D$ and $D \bigcap h(D)$  and two translations arcs, these four arcs having disjoint interiors.}

\medskip

{\it Proof}. Because of Jordan's theorem, $h^{-1}(D) \bigcap D$ and $D \bigcap h(D)$ are not linked on $\partial D$ (that is, are contained in disjoint intervals of $\partial D$). Let $\lambda$ be the smallest closed sub-arc (eventually reduced to a point) of $\partial D$ containing $h^{-1}(D) \bigcap D$ and disjoint from $h(D)$. According to the preceding Lemma and Jordan's Theorem, $h(D) \bigcap D$ is contained either in int$\lambda$ or $\partial D \setminus \lambda$. In the first case, the frontier of the unbounded component of ${\bf R}^2 \setminus (h^{-1}(D) \bigcup D)$ is a Jordan curve send into the disc it bounds by $h$, which gives a fixed point for $h$: a contradiction.
Therefore, the second case holds and the arc $\mu$, smallest closed sub-arc (eventually reduced to a point) of $\partial D$ containing $D \bigcap h(D)$, is disjoint from $\lambda$. The end points of $\mu$ and $h(\lambda)$ coincide and, since $h$ preserves orientation, the two sub-arcs of $\partial D$ which are the components of $\partial D \setminus (\hbox{int}\lambda \bigcup \hbox{int}\mu )$ are translation arcs.

\bigskip

{\bf 2.3 Definition}. A translation arc as in the preceding lemma, oriented as said before, will be said {\bf high} if $D$ is on the right of the translation arc and {\bf low} if $D$ is on the left of the translation arc.

\bigskip

\noindent{CONVENTION. From now on, unless said otherwise, all discs will be {\bf euclidean} discs.}

\bigskip

{\bf 2.4 Definition}. Let us consider two critical discs $D'$ and $D$ which satisfy:

1.$\partial D \bigcap D'$ is a non degenerated sub-arc of a translation arc $\alpha $ of $\partial D$.

2. $h^{-1}(D) \bigcap \hbox{int}D' = \emptyset = \hbox{int}D' \bigcap h(D)$.

$D'$ will be called an {\bf extension} of $D$ (high or low according to the type of the arc $\alpha$).

\bigskip

{\bf 2.5 Remarks}.

1. $D'$ is a high extension of $D$ if and only if $D$ is a low extension of $D'$. Indeed, it is enough to check that the low translation arc of $D'$ contains $\partial D' \bigcap D$. But, if this were not true, int$D$ would contain a point of $D' \bigcap h(D')$ or of $D' \bigcap h^{-1}(D')$ in contradiction to point 2 of the definition.

2. $D \bigcup D'$ is still a critical  topological disc when $D'$ is a high extension of $D$.

\bigskip

{\bf 2.6 Lemma}. {\it Let $D$ and $D'$ be two critical discs and $\alpha _-$, $\alpha _+$ the two arcs of translations of $\partial D$. Suppose that int$D' \bigcap \partial D$ is not empty, meets $\alpha _+$ and that $D' \bigcap (h^{-1}(\alpha _+) \bigcup h(\alpha_+)) = \emptyset$.

Then, $D \bigcup D'$ is critical. Therefore, $D'$ is a high extension of $D$.}

\bigskip

{\bf 2.7 Remark}. Of course, we have an analogous Lemma obtained by exchanging $\alpha _+$ and $\alpha _-$ and ``high'' and ``low''.

\medskip

{\it Proof}. We will first prove:

\medskip

{\bf 2.8 Assertion}. {\it Let $D$ be any critical or free topological disc and $\alpha$ any translation arc such that int$D \bigcap \alpha \not = \emptyset$ and $D \bigcap (h^{-1}(\alpha) \bigcup h(\alpha)) = \emptyset$, then $D \bigcap \bar l = D \bigcap \alpha$ where $l$ is the trajectory generated by $\alpha$.}

\medskip

Given this assertion and the fact that $l$ separates the connected component of ${\bf R}^2 \setminus (\bar l \setminus l)$ containing $l$ into two components invariant under $h$, one concludes easily the proof of the Lemma.

\medskip

{\it Proof of the Assertion}. Suppose first that $D \bigcap h^n (\alpha) \not = \emptyset$ for some $n \in {\bf Z}$, $\vert n \vert >1$. We can find an arc $\gamma$ in int$D$, except perhaps for one end point, joining a point $x \in \hbox{int}\alpha \bigcap \hbox{int}D$ to a point $h^n(y)$ for some $y \in \alpha$. Since int$D$ and the subarc $[x, y]$ of $\alpha$ are free, $\gamma$ is free also. Then, small enough open topological discs, neighborhoods of $[h(x), h(y)]$ and $\gamma$ respectively,  would satisfy the hypothesis of Franks' Lemma giving rise to a fixed point for $h$ and a contradiction.

If $D \bigcap \bar l \not = \emptyset$, we can suppose that $D \bigcap l = \emptyset$  from what precedes and we can find a free arc $\gamma$ inside int$D$, except perhaps for one end point, joining a point $x \in \hbox{int}\alpha \bigcap \hbox{int}D$ to a point $z \in \bar l \setminus l$. We choose then a small free disc $D'$ around $z$ such that $D' \bigcap (h^{-1}(\gamma) \bigcup h(\gamma)) = \emptyset$. Then $D' \bigcup \gamma$ is free and a free topological disc, small enough neighborhood of $D' \bigcup \gamma$, verifies the same hypothesis as $D$ in the first paragraph.\cqfd

\bigskip

{\bf 2.9 Lemma}. {\it If $(D_i)_{i \in I}$ (where $I$ is a finite set $\{ 0,1, \ldots , n\}$ or ${\bf N}$ or ${\bf Z}$) is a sequence of critical discs such that $D_{i+1}$ be a high extension of $D_i$, then $\bigcup _{i \in I}D_i$ is a critical disc. Furthermore, $(D_i \setminus D_{i-1}) \bigcap (D_j \setminus D_{j-1}) = \emptyset$ if $i \not = j$.}

\medskip

{\it Proof}. Let $l_i$ be the translation line generated by the high translation arc of $D_i$. By Assertion 2.8, the high translation arc of $D_{i+1}$ does not meet $l_i$ and we define $\Sigma _i$ as the component of ${\bf R}^2 \setminus (\bar l_i \setminus l_i)$ which contains the high translation arc of $D_{i+1}$. It follows from the preceding Lemma that $\Sigma _i \supset \Sigma _{i+1}$. Using the invariance of these $\Sigma_i$ under $h$ one gets the Lemma by an easy induction (and, in case $I = {\bf Z}$, a second induction going from $i$ to $i-1$ and using the low translation arc of $D_i$).

\bigskip

To each point $x$ in the plane, we now associate the radius $r_x$ of the euclidean critical disc with center at $x$.

\medskip

{\bf 2.10 Lemma}. {\it The map $x \mapsto r_x$ is continuous.}

\medskip

{\it Proof}. This follows from the fact that, if $d(x, y)< \epsilon$, then $B(y, r_x -\epsilon) \subset B(x, r_x ) \subset B(y, r_x + \epsilon)$.

\bigskip

{\bf 2.11 Lemma}. {\it Let $D$ be a critical disc and $\alpha$ be one of the two translation arcs  of $\partial D$. There exists a  critical disc $D'$ with center on $\alpha$ such that $D'$ is an extension (high or low according to the choice of $\alpha$) of $D$.}

\medskip

{\it Proof}. Let $C_+$ (resp. $C_-$) be the set of points $x$ in $\alpha$ such that the critical disc with center at $x$ meets $h(\hbox{int}D)$ (resp. $h^{-1}(\hbox{int}D)$). To prove the Lemma it is enough to prove that $C_+ \bigcup C_-$ is not all of $\alpha$. This will follow from the connectivity of $\alpha$ and the following three points.

$\bullet$ $C_+ \bigcap C_- = \emptyset$: Let $D'$ be an euclidean critical disc in $C_+ \bigcap C_-$, that is, meeting $h(\hbox{int}D)$ and $h^{-1}(\hbox{int}D)$. Then there exists a topological disc $\Delta \subset \hbox{int}D' \setminus \hbox{int}D$ meeting both $h(\hbox{int}D)$ and $h^{-1}(\hbox{int}D)$ which gives with int$D$ two discs satisfying the hypothesis of Franks' Lemma and therefore a contradiction.

$\bullet$ $C_+$ and $C_-$ are not empty: it is enough to consider discs centered on the end points of $\alpha$.

$\bullet$ $C_+$ and $C_-$ are open: this follows from the continuity of the map $x \mapsto r_x$.

\bigskip

Under the hypothesis and notations of the preceding Lemma, since $D'$ is critical, one can repeat the construction using this time the translation arc on $\partial D'$ of the same type than $\alpha$. One gets in this way an infinite sequence $D_0 = D, D_1 = D', D_2, \ldots , D_n, \ldots$ of critical discs such that by Lemma 2.9, for every $n$, $\bigcup _{i=0}^{i=n}D_i$ is a critical disc and $(D_i \setminus D_{i-1}) \bigcap (D_j \setminus D_{j-1}) = \emptyset$ if $i \not = j$.

\bigskip

{\bf 2.12 Lemma}. {\it In this situation, the sequence $\{D_n \}$ converges to infinity, that is, given any compactum $K \subset {\bf R}^2$, there exists $N$ such that the connected set $\bigcup _{i\geq N} D_i$ does not meet $K$.}

\medskip

{\it Proof}. Since the $D_n$ are critical, the radius of those meeting $K$ is bounded below as well as their area. By the Lemma all $D_n \setminus D_{n-1}$ are disjoint and, by euclidean geometry, have area greater than half the area of $D_n$. Therefore only a finite number of discs $D_n$ can meet $K$.

\bigskip

It is now easy to prove the Brouwer plane translation theorem.

\medskip

{\bf 2.13 Definition}. A {\bf Brouwer line} $L$ for the Brouwer homeomorphism $h$ is the image of a proper embedding of ${\bf R}$ into ${\bf R}^2$ such that

1. $L$ is free under $h$.

2. $L$ separates $h(L)$ and $h^{-1}(L)$.

\bigskip

{\bf 2.14 Plane Translation Theorem (Brouwer)}. {\it Let $h$ be a Brouwer homeomorphism. Then every point in ${\bf R}^2$ belongs to a Brouwer line.}

{\it Proof}. Let $x$ be the point in consideration and $D_0$ the critical disc centered at $x$. One constructs as above a sequence $D_0, D_1, \ldots , D_n, \ldots $ of critical discs converging to infinity, and also, starting from the other translation arc of $D_0$, another sequence  $D_0, D_{-1}, \ldots , D_{-n}, \ldots $ of critical discs converging to infinity. Note that $\bigcup _{k=i}^{k=j}D_k$ is critical for $i<j$ where $i, j \in {\bf Z}$ since one can see that finite sequence as obtained by the preceding construction from $D_i$ and therefore $\bigcup _{i \in {\bf Z}}\hbox{int}D_i$ is free. It is then easy to construct inside $\bigcup _{i \in {\bf Z}}\hbox{int}D_i$ a proper embedding $L$ of ${\bf R}$ joining linearly the centers of the $D_i$'s. It has to be free and also to satisfy property 2 of a Brouwer line since one can find in $D_0$ an simple arc joining $h^{-1}(D_0)$ to $h(D_0)$ crossing $L$ only at $x$ and if $L$ where not separating $h^{-1}(L)$ and $h(L)$, one could complete this arc into a simple closed curve crossing $L$ only at $x$ in contradiction to the properness of $L$.

\vskip 2cm

\noindent{\bf 3. Further results on euclidean critical discs}.

\bigskip

We still fix a Brouwer homeomorphism $h$.

\medskip

{\bf 3.1 Definition}. Let $D$ be a critical disc. We know that its boundary is the union of two translation arcs $\alpha_+$ and $\alpha_-$ (high and low respectively) and of two disjoints arcs, $\lambda$ and $\mu$, smallest sub-arcs of $\partial D$ containing $h^{-1}(D) \bigcap D$ and $D \bigcap h(D)$ respectively. Choose $\epsilon$ so that $0< \epsilon \leq {\displaystyle {1 \over 8}} \hbox{inf}(l(\alpha_+), l(\alpha_-), l(\lambda), l(\mu))$ where $l(.)$ is the length of the corresponding arc and where we omit $l(\lambda )$ and $l(\mu )$ if $\lambda$ (or equivalently $\mu$) is degenerated. We have $\epsilon < {\displaystyle {\pi \over 8}}r$ if $r$ is the radius of $D$.

Then $D$ without the discs of radius $\epsilon$ centered at the end points of $\alpha_-$ is a topological disc $D_+$ and $D$ without the discs of radius $\epsilon$ centered at the end points of $\alpha_+$ is also a topological disc $D_-$. If $\lambda$ (or $\mu$) is not degenerated $D_+$ and $D_-$ are topological critical discs and othervise $D_+ = D_-$ is a free topological disc. In any case, area$(D_{\pm}) > {\displaystyle{{2 \pi} \over 3}}r^2$.

So, to each point $x$ of the plane we may associate three critical topological discs: the euclidean critical disc $D$ centered at $x$ and, if $\lambda$ and $\mu$ are non degenerated, two {\bf derived} topological discs $D_+$ and $D_-$. We will also say that $D$ is derived from $D_+$ or $D_-$. These three discs will be called {\bf generalized critical discs}. We will also say that $D$ is derived from $D_+$ or $D_-$. The centers, radii and translation arcs of the derived discs are those of $D$ by definition. We will extend this vocabulary to the case where $\lambda$ and $\mu$ are degenerated: in that case $D_+ = D_-$ is free and not critical but, having its area bounded below as said above, it will not cause any problem below.

The euclidean critical disc $D$ will be said {\bf the underlying disc} of $D_+$ and $D_-$ and also $D$ on occasion.

\bigskip

We will begin to improve the extension Lemma as follows.

\bigskip

{\bf 3.2 Lemma}. {\it Let $D$ be a critical disc and $\alpha$ one of the two translation arcs of $\partial D$. Then there exists a generalised critical disc $D'$ such that

1. The underlying disc of $D'$ is an extension of $D$.

2. $D' \bigcap \partial D = D' \bigcap \hbox{int}\alpha$

3. $h^{-1}(D) \bigcap D' = \emptyset = D' \bigcap h(D)$}

\bigskip

{\bf 3.3 Definition}. Given two critical discs $D_1$ and $D_2$ admitting derived discs $D$ and $D'$ which satisfy conditions 1, 2 and 3 above we will say that $D_2$ is a {\bf strict extension} of $D_1$ (high or low according to the type of $\alpha$).

\bigskip

{\bf 3.4 Remarks}. 1. $D \bigcup D'$ is a topological critical disc (or free if $D$ and $D'$ are free).

2. If $D_2$ is a strict high (resp. low) extension of $D_1$, then $D_1$ is a strict low (resp. high) extension of $D_2$.

\bigskip

{\it Proof}. Let $A_+$ (resp. $A_-$) be the subset of $\alpha$ of points $x \in \alpha$ such that the critical disc with center at $x$ meets $h(D)$ (resp. $h^{-1}(D)$). Then

$\bullet$ $A_+$ and $A_-$ are not empty: it is enough to look at critical discs centered at the end points of $\alpha$.

$\bullet$ $A_+$ and $A_-$ are closed: this follows from the continuity of the map $x \mapsto r_x$.

$\bullet$ If $A_+ \bigcap A_- = \emptyset$, the connectivity of $\alpha$ implies that an honest critical disc satisfies the Lemma.

$\bullet$ If $A_+ \bigcap A_- \not = \emptyset$, then this intersection is reduced to only one disc: the one centered at the midpoint of $\alpha$ and with the end points of $\alpha$ on its boundary. 

Indeed, if there were a disc $D'$ meeting $h^{-1}(D)$, $h(D)$ and at most one end point of $\alpha$, one could find a free arc $\delta$ contained in $D' \setminus D$ except for its endpoints $a \in h^{-1}(D)$ and $b \in h(D)$. One could also find a free arc $\gamma$ in $D$, disjoint from $\delta$ and joining $h(a)$ to $h^{-1}(b)$ and therefore, by fattening, a contradiction with Franks'Lemma.

Note that the same argument shows that a critical disc $D$ which contains the two end points of $\alpha$ cannot meet $h^{-1}(D)$ or $h(D)$ except in these end points.

Then that one disc $D'$ admits a derived disc ($D'_+$ or $D'_-$ according to the type of $\alpha$) which satisfies the first case of the Lemma.

\bigskip

{\bf 3.5 Lemma}. {\it Analogs of Lemmas 2.9 and 2.12 are now equally true with extension replaced by strict extension and with the same proofs (for the second Lemma we use the fact that the free discs eventually met have an area bounded below as said above).}

\bigskip

It is only to get the following Lemma that we had to consider strict extensions and derived discs. But of course this Lemma will be crucial in section $4$.

\bigskip

{\bf 3.6 Lemma}. {\it Let $K$ be a compactum in the plane. There exists a finite number $n$ of pairs of generalised discs $(D_i^1, D_i^2), 1 \leq i \leq n$, such that every critical disc centered on $K$ admits one of these pairs as a couple of strict extensions: high and low.}

\medskip

{\it Proof}. If $D$ is an euclidean critical disc centered at $x$ and $(D^1, D^2)$ a pair of euclidean generalised critical discs such that $D^1$ is a strict high extension of $D$ and $D^2$ a strict low extension of $D$, then there exists a neighborhood $V_x$ of $x$ such that the euclidean critical discs centered on $V_x$ all admit $D^1$ and $D^2$ as strict high and low extensions (since $x \mapsto r_x$ is continuous and the arcs $\partial D \bigcap D^1$ and $\partial D \bigcap D^2$ being free, neighboring arcs are still free). A finite sub-covering of the covering of $K$ by the $V_x$'s leads to the sought after pairs.

\vskip 2cm

\noindent{\bf 4. Homeomorphisms of the torus and the open annulus}.

\bigskip

We continue to call $h$ a Brouwer homeomorphism of the plane.

\bigskip

We will denote by $s$ and $t$ the applications $(x, y) \mapsto (x, y+1)$ and $(x, y) \mapsto (x+1, y)$ of ${\bf R}^2$. Also, $\pi$ will be the projection ${\bf R}^2 \rightarrow {\bf R}/{\bf Z} \times {\bf R}$.

\bigskip

{\bf 4.1 Definition}. If $L$ is a Brouwer line for $h$, we denote by $D(L)$ the component of ${\bf R}^2 \setminus L$ which contains $h(L)$ and by $G(L)$ the one which contains $h^{-1}(L)$.

\bigskip

{\bf 4.2 Lemma}. {\it Let $(L_n)_{n\in {\bf Z}}$ be a locally finite family of Brouwer lines for $h$ such that $\overline {\bigcup _n D(L_n)}$ is connected and ${\bf R}^2 \setminus \overline {\bigcup _n D(L_n)}$ is a connected, unbounded open set $R$. Then, the frontier of $R$ is a Brouwer line for $h$.}

\medskip

{\it Proof}. Given the local finiteness of the family, $\overline {\bigcup _n D(L_n)} = \bigcup _n \overline {D(L_n)}$ and a classical result of K\'er\'ekjart\'o [K1, page 87] implies that $R$ is a submanifold of the plane. ${\bf R}^2$ being unicoherent and  $\overline {\bigcup _n D(L_n)}$ connected, Fr$R$ is connected. But the frontier of $R$ is free, because if $x \in \hbox{Fr}R$, there is $n$ such that $x \in L_n$ and therefore, $h(x) \in D(L_n) \subset {\bf R}^2 \setminus \overline {R}$. Since $h^{-1}(R) \subset R$, it cannot be a circle and so it is a proper free line. It is a Brouwer line since $h^{-1}(R) \subset R$ and  $h(R) \subset {\bf R}^2 \setminus \overline {R}$.

\bigskip

{\bf 4.3 Lemma}. {\it Suppose $h$ commutes with $t$ and let be given a Brouwer line $L$ for $h$ which projects onto and properly on the $y$-axis. Then there exists a Brouwer line $L'$ for $h$ which projects onto and properly on the $y$-axis and free for $t$.}

{\it Furthermore, if $h$ commutes also with $s$ and if $L$ is invariant under $s$, one can find $L'$ also invariant under $s$.}

\bigskip

{\it Proof}. Since $L$ projects down properly onto the $y$-axis, there exists $N>0$ such that $G(L)$ or $D(L)$ contains $[N, +\infty ) \times \{ 0 \}$. Exchanging $h$ and $h^{-1}$ we can suppose that $D(L)$ contains a half line $[N, +\infty ) \times \{ 0\}$ and consider the unbounded connected component $R$ of ${\bf R}^2 \setminus \bigcup_{k \geq 0}t^k(D(L))$ and its frontier $L'$.

For every $A>0$ there is $n(A)$ such that for $n \geq n(A)$, $t^n(L \bigcap {\bf R} \times [-A, A]) \bigcap {\bf R} \times [-A, A] = \emptyset$, therefore the family $(t^k(L))_{k \geq 0}$ is locally finite and $L'$ projects onto and properly on the $y$-axis. Furthermore, by Lemma 4.2 , $L'$ is a Brouwer line and it satisfies $t(L') \subset {\overline D(L')}$ since if $x \in L'$, there exists $n$ such that $x \in t^n(L)$ and so $t(x) \in t^{n+1}({\overline D(L)} \subset {\bf R}^2 \setminus R = {\overline D(L')}$.

Notice that, since $s$ and $t$ commute, if $s(L)=L$ then also $s(L')=L'$.

Therefore to conclude that, in this second case, exists a Brouwer line going down in the annulus on a free line joining the two ends, it is enough to apply the following Sublemma to $L'$.

\bigskip

{\bf 4.4 Sublemma}. {\it Let $L$ be a Brouwer line for $h$ with a proper onto projection on the $y$-axis such that $t(L) \subset {\overline D(L)}$. Then there exists a Brouwer line $L'$ with a proper onto projection on the $y$-axis such that $t(L') \subset D(L')$.}

{\it Furthermore, if $h$ commutes with $s$ and $s(L)=L$, then one can choose $L'$ such that $s(L')=L'$.}

\medskip

{\it Proof}. Let $\delta : L \rightarrow ]0, +\infty[$ be a function such that the $\delta$-neighborhood $W$ of $L$ does not meet $h^{-1}(W) \bigcup h(W)$ and define $L_0 = \emptyset$ and, for all $k>0$, $L_k = L \bigcap ({\bf R} \times [-k, k])$. Since $L_k$ is compact, there exists $ n = n(k)$ such that $L_k \bigcap t^n(L_k) = \emptyset$. Let us suppose that $L_{k-1} \bigcap t(L_{k-1}) = \emptyset$ and let us show that there exist a Brouwer line $L'$ with a proper onto projection on the $y$-axis such that $L'_{k-1}=L_k$ and $L'_k \bigcap t^m(L'_k) = \emptyset$ for some $m, 0<m<n$, if $n>1$. (This is enough to finish the proof of the Sublemma).

\medskip

Define $X_i = L_k \bigcap \ldots \bigcap t^i(L_k)$ if $i\geq 0$. By hypothesis, $t(X_{n-1}) \bigcap L_k = X_n = \emptyset$ and $X_{n-1}$ is a compact of $L_k$ disjoint from $L_{k-1}$ as well as the $\delta$-neighborhood of $X_{n-1}$ (by diminishing $\delta$ near $X_{n-1}$ if necessary). The $\delta$-neighborhood of $X_{n-1}$ in $L$ is a finite union of open intervals in $L_k \setminus L_{k-1}$. Let us replace these intervals by arcs with the same end points but included in the intersection of $G(L)$ and of the $\delta$-neighborhood of $L$. We thus get a Brouwer line $L'$ which projects onto and properly on the $y$-axis and such that $L'_{k-1} = L_{k-1}$. We are left to check that $X'_{n-1} = \emptyset$.

\medskip

To do that, notice that $t(L_k) \bigcap L'_k \subset X_1 \setminus X_{n-1}$ and that $t(L'_k)$ and $t(L_k)$ differ only in the $\delta$-neigborhood of $t(X_{n-1})$. For $\delta$ small enough, this neighborhood does not meet $L_k$ or $L'_k$ and therefore $L'_k \bigcap t(L'_k) = L'_k \bigcap t(L_k) \subset X_1 = L_k \bigcap t (L_k)$. We deduce that $t^i(L'_k) \bigcap t^{i+1}(L'_k) \subset t^i(L_k) \bigcap t^{i+1}(L_k)$ for $i \geq 0$ and $X'_{n-1} \subset X_{n-1}$ and conclude that $X'_{n-1} = L'_k \bigcap t(L'_k) \bigcap X'_{n-1} \subset L'_k \bigcap t(L'_k) \bigcap X_{n-1} = L'_k \bigcap t(L_k) \bigcap X_{n-1} \subset (X_1 \setminus X_{n-1}) \bigcap X_{n-1} = \emptyset$.

\medskip

In case $s(L)=L$, one can choose the function $\delta$ such that $\delta s = \delta$ so that $s(W)=W$. We forget all $L_k, k \geq 0$, consider instead only $L$ and note that since $s(L)=L$ there exists $n>0$ such that $L \bigcap t^n(L)= \emptyset $. We define analogously $X_i = L \bigcap \ldots \bigcap t^i(L)$ for $i\geq 0$. Since $s$ and $t$ commute, each $X_i , i\geq 0,$ satisfies $s(X_i)=X_i$, therefore we can choose the family of arcs changing $L$ to $L'$ so that $s(L')=L'$. \cqfd

\bigskip

{\bf 4.5 Theorem}. {\it Let $H$ be a homeomorphism of the torus ${\bf T}^2$ preserving orientation which admits a lift $h$ to ${\bf R}^2$ without fixed point. Then there exists an essential simple closed curve on ${\bf T}^2$ which lifts to a Brouwer line of $h$.}

\bigskip

{\it Proof}. Applying Lemma 3.6, we get a finite number of generalised critical discs $\{ D_1, \ldots , D_r \}$ such that every critical disc centered on $[0, 1] \times [0, 1]$ admits one of the $D_i$'s, $1 \leq i \leq r$ as high strict extension. Since $h$ commutes with $s$ and $t$, we have in fact a  finite set of euclidean generalised critical discs  $\{ D_1, \ldots , D_r \}$ such that for every critical disc $D$ of ${\bf R}^2$ there exists $n, m \in {\bf Z}$ and $i$ between $1$ and $r$ such that $s^nt^m(D_i)$ is a high strict extension of $D$.

We choose then a critical disc $D$. To this disc we associate a disc $D'$ from $\{ s^nt^m(D_i) \vert n, m \in {\bf Z}, 1 \leq i \leq r \}$ as just said. Since the underlying disc to $D'$ is critical, there exists also $D''$ in the same set which is a high extension of $D'$. In that way, we construct an infinite sequence $D', D'', \ldots , D^{(n)}, \ldots $ of generalised critical discs such that $\bigcup _{i=1}^{i=n}D^{(i)}$ is a critical topological disc for every $n$ (according to Lemma 3.5). Now, the finiteness of  $\{ D_1, \ldots , D_r \}$ implies that there exists $n, m \in {\bf Z}$ and $j > i \geq 1$ such that $D^{(j)} = s^nt^m(D^{(i)})$. We note that $C = \bigcup _{p\in {\bf Z}}(s^nt^m)^p(D^{(i)} \bigcup \ldots \bigcup D^{(j-1)})$ is a critical band since $(s^nt^m)^{p+1}(D^{(i)}) = (s^nt^m)^p(D^{(j)})$ is a high extension of $(s^nt^m)^p(D^{(j-1)})$. Of course, $s^nt^m(C) = C$ and we can find a Brouwer line $L$ inside $C$ such that $s^nt^m(L)=L$. If $n$ and $m$ are relatively prime (or if $n=0$ and $m=1$ or if $m=0$ and $n=1$) we can conjugate $h$ by an element of $SL(2, {\bf Z})$ and suppose that $m=0, n=1$, that is, that $s(L)=L$. Lemma 4.3 gives then a Brouwer line $L'$ such that $s(L') =L'$ and $t(L') \bigcap L' = \emptyset$ which therefore goes down to a simple closed curve in ${\bf R} \times S^1$ free under the map induced by $t$ and then to a simple closed curve in ${\bf T}^2$. In case $d=\hbox{gcd}(m, n)>1$, we  apply  Lemma 4.2 to the Brouwer lines $L, T(L), \ldots , T^{d-1}(L)$ where $T = s^{n \over d}t^{m \over d}$ is a translation of ${\bf R}^2$ such that $T^d(C)=C$ since $s$ and $t$ commute (if $n=0$ or $m=0$ we let $d = \vert n \vert + \vert m \vert$). The Brouwer line given by this Lemma (Fr$R$ there) is invariant by the translation $T$ and we are back to the preceding case. \cqfd

\bigskip

{\bf 4.6 Lemma}. {\it There is no euclidean free or critical half plane if $h$ is a lift to ${\bf R}^2$ of a homeomorphism $H$ of $S^1 \times {\bf R} = {\bf R}/{\bf Z} \times {\bf R}$ isotopic to the identity.}

\medskip

{\it Proof}. Suppose first that the boundary of the half plane $D$ is not parallel to the $x$-axis. Then, given $z \in {\bf R}^2$, there exists $n \in {\bf Z}$ such that $t^n(z)$ and $t^n(h(z))$ are in int$D$, but $t^n(h(z)) = h(t^n(z))$ and so $D$ is not free or critical.

If the boundary of the free or critical half plane $D$ is parallel to the $x$-axis, $H(\pi (\hbox{int}D) \bigcap \pi (\hbox{int}D) = \pi (h(\hbox{int}D) \bigcap \hbox{int}D) = \emptyset$ (the first equality follows from the $t$-invariance of $D$). But $\pi (\hbox{int}D)$ is a neighborhood of an end of $S^1 \times {\bf R}$ and $H$ must preserve the ends: contradiction.

\bigskip

{\bf 4.7 Theorem}. {\it Let $H$ be a homeomorphism of the open annulus $S^1 \times {\bf R}$ isotopic to the identity and with a fixed point free lift $h$ to ${\bf R}^2$. Then,

there exists an essential simple close curve in the annulus free under $H$ (and therefore lifting to a Brouwer line for $h$ in ${\bf R}^2$)

or there exist a properly embedded line in the annulus going from one end to the other which lifts to a Brouwer line for $h$.}

\medskip

{\it Proof}. Applying Lemma 3.6 to $[0, 1] \times [n, n+1] \subset {\bf R}^2$, we get a finite set of generalised critical discs $\{ D_1^n, \ldots , D_{r_n}^n \}$ such that every critical disc centered on $[0, 1] \times [n, n+1]$ admits one of the $D_i^n, 1 \leq i \leq r_n$ as high strict extension and another one as low strict extension. Since $h$ and $t$ commute, for every critical disc $D$ centered on the band ${\bf R} \times [n, n+1]$ there exist a $t^k(D^n_i), 1 \leq i \leq r_n , k\in {\bf Z}$, which is a high strict extension of $D$ and another one  which is a low strict extension. We can suppose, by adding, if necessary, a finite number of critical discs to each family $\{ D_1^n, \ldots , D_{r_n}^n \}$, that ${\cal D} = \bigcup _{n \in {\bf Z}}\{ D_1^n, \ldots , D_{r_n}^n \}$ covers $[0, 1] \times {\bf R}$ and so, that $\pi ({\cal D})$ covers $A$. 

\bigskip

{\bf 4.8 Fact}. {\it Only a finite number of discs of the family ${\cal D}$ can meet a given band $B_n = {\bf R} \times [n, n+1]$.}

\medskip

{\it Proof}. First remark that critical discs centered on $[0, 1] \times {\bf R}$ with radius going to infinity have to go to infinity. Because, on the contrary, since the centers of these discs have to go to infinity, we would get a free or critical euclidean half plane contrarily to Lemma 4.6. 

Suppose then that an infinite number of discs $D'_k$ in ${\cal D}$ meet the band $B_n$. By definition $D'_k$ meets also a critical disc $D_k$ centered on $[0, 1] \times [n_k, n_k +1]$ where $n_k \rightarrow \pm \infty$. If the radius of $D_k$ is bounded, the radius of $D'_k$ has to go to infinity, but if the radius of $D_k$ is not bounded, we can suppose (passing to a subsequence) that it goes to infinity and, using the preceding remark, we conclude that the radius of $D'_k$ has still to go to infinity. Since translates of critical discs are still critical, we can get a sequence of critical discs which have their centers on $[0, 1] \times {\bf R}$, meet $B_n$ and have a radius going to infinity: this implies the existence of a free or critical euclidean half plane contrarily to Lemma 4.6. \cqfd

\bigskip

We then choose a critical disc $D$. There is a disc $D'\in  \bigcup _{k,n \in {\bf Z}}t^k(\{ D_1^n, \ldots , D_{r_n}^n \})$ which is a high strict extension of $D$. Then there exists $D''$ in the same set which is a high strict extension of $D'$ and so on ad infinitum. One gets an infinite sequence $D', D'', \ldots , D^{(n)}, \ldots $ of generalised critical discs such that $D \bigcup (\bigcup _{i>0}D^{(i)})$ is a proper critical half band (see Lemma 3.5). The same reasoning with low extensions leads to a proper critical band $D_{\infty} = (\bigcup _{i<0}D^{(i)}) \bigcup  (\bigcup _{i>0}D^{(i)})$ made of generalised critical discs.

We will consider three cases covering all possibilities for the relative positions of the band $D_{\infty}$ and the bands $B_n = {\bf R} \times [n, n+1]$.

\medskip

{\bf First case}. There is an infinite number of the discs $D^{(n)}$ which meet a fixed band $B_k$.

As a consequence of Fact 4.8 above, there exist $k \in {\bf Z}$ and $j>i \geq 1$ such that $D^{(j)}=t^k(D^{(i)})$ so that $t^{k(l+1)}(D^{(i)})=t^{kl}(D^{(j)})$ is a high extension of $t^{kl}(D^{(j-1)})$ and consequentely $C = \bigcup _{\l \in {\bf Z}}t^{kl}(D^{(i)}) \bigcup \ldots \bigcup D^{(j-1)})$ is a critical band. Of course $t^k(C) =C$ and $C$ contains a Brouwer line $L$ such that $t^k(L)=L$. Applying Lemma 4.2, we see that there exist an essential simple closed curve as in the Theorem.

\medskip

{\bf Second case}. For all $n\in {\bf Z}$, the band $D_{\infty}$ contains a non zero but finite number of $D^{(i)}$ which meet the band $B_n$.

Therefore there exists a Brouwer line $L$ which projects onto and properly on the $y$-axis and we conclude using Lemma 4.3.
\medskip

{\bf Third case}. There exists $n$ such that $D_{\infty}$ does not contain any disc meeting the band $B_n$. In other words, $D_{\infty}$ is contained in an euclidean half plane with horizontal boundary.

\medskip

We can suppose that, for every $n \in {\bf Z}$, $\bigcup _{i<0}D^{(i)}$ et $\bigcup _{i>0} D^{(i)}$ contain only a finite number of discs meeting $B_n$ for otherwise we are in the first case. We distinguish two sub cases:

\medskip

Either $t(D_{\infty}) \bigcap D_{\infty} \not = \emptyset$. Consider then a Brouwer line $L \subset D_{\infty}$ such that $t(L) \bigcap L \not = \emptyset$. An application of Lemma 4.2 to the family $(t^n(L))_{n \in {\bf Z}}$ of Brouwer lines leads to an essential simple closed curve as required in the Theorem.

\medskip 

Or $t(D_{\infty}) \bigcap D_{\infty} = \emptyset$. In that case, by Lemma 2.1, $t^n(D_{\infty}) \bigcap D_{\infty} = \emptyset$ for all $n \not = 0$ and we can find a Brouwer line $L$ such that $t^n (L) \bigcap L = \emptyset$ for all $n \not =0$. Such a line goes down in the annulus onto a free proper line going from one end to itself.

\medskip

To conclude the proof of the Theorem it is therefore enough to show the following:

\medskip

{\bf 4.9 Affirmation}. {\it In the preceding situation, suppose that any band $D_{\infty}$ constructed as above starting from any critical disc of the family ${\cal D} = \bigcup_n \{ D_1^n, \ldots , D^n_{r_n} \}$ satisfy $t^n(D_{\infty}) \bigcap D_{\infty} = \emptyset$ for all $n \not = 0$. Then

there exists in the annulus an essential simple closed curve free under $H$

or there exists a line in the annulus joining the two ends of the annulus and lifting to a Brouwer line for $h$.}

\medskip

{\it Proof}. We denote by $N$ and $S$ the ends of the annulus $A = {\bf R}/{\bf Z} \times {\bf R}$. Choose for each $D \in {\cal D}$ a band $D_{\infty}$ constructed as above. Given the hypothesis there are four types of generalized critical discs in  ${\cal D} = \bigcup_n \{ D_1^n, \ldots , D^n_{r_n} \}$: either the band $D_{\infty}$ associated to such a disc verifies that $\pi (D_{\infty})$ goes from $N$ to $N$ and, according to the separation or not of $H(D_{\infty})$ and $S$ by $D_{\infty}$, we get the types $\rightarrow \! N$ and $\leftarrow \! N$, or $\pi (D_{\infty})$ goes from $S$ to $S$ and we get the types $\leftarrow \! S$ and $\rightarrow \! S$.

\medskip

$\bullet$ Let us suppose now that there exist discs $D \in \cal D$ of type $\rightarrow \! S$ or $\leftarrow \! S$ such that the $\pi (D)$'s meet arbitrarly small neighborhoods of $N$. We choose then a sequence $(D_n)_{n \geq 1}$ of such discs such that $\pi (D_n)$ meets ${\bf R}/{\bf Z} \times [n, +\infty [$. To each $D_n$ is associated a properly embedded band $\pi (D^n_ {\infty })$ in the annulus whose two ends converge towards the end $S$ of the annulus $A$. Let $F_k$ the union over $n$ of the set of connected components of $\pi (D^n_{\infty} \bigcap ({\bf R}/{\bf Z} \times [-k, k]))$ which join ${\bf R}/{\bf Z} \times \{ -k \}$ to  ${\bf R}/{\bf Z} \times \{ k \}$. The $F_k$ are finite sets by Fact 4.6 and every element of $F_{k+1}$ restrict to one of $F_k$. By the set theoretical Lemma below, there exists a band of $A$ going from $N$ to $S$, made from a ${\bf Z}$ indexed sequence of elements of $\pi ({\cal D})$, every element of this sequence being a high extension of the preceding one. Therefore, this band is critical for $H$ and we get a properly embedded free line in the annulus $A$ joining $N$ and $S$ as desired.

Of course, we get the same conclusion if we suppose that there exist discs $D \in \cal D$ of type $\rightarrow \! N$ or $\leftarrow \! N$ such that the $\pi (D)$ meet arbitrarly small neighborhoods of $S$.

\medskip

$\bullet$ Let us suppose that the two cases considered in the preceding point do not show up and particularly that all discs $\pi (D) \in \pi (\cal D)$ with $D$ of type $\rightarrow \! N$ avoid a neighbohood $O_S$ of $S$ and that all those of type $\rightarrow \! S$ avoid a neighborhood $O_N$ of $N$.

Let us first remark that, if $D \in \cal D$ is of type $\rightarrow \! S \ (\hbox{resp.} \rightarrow \! N)$, every point of $\pi (D)$ converges towards $S$ (resp. $N$) under the action of $H$. Indeed, if $D_{\infty}$ is the critical band associated to $D$, then all $h^n(D)$ are separated from $N$ by $D_{\infty}$ and the set of points in the plane separated from $N$ by $D_{\infty}$ and with $y$-coordinate bounded below is bounded. But all points of the plane are wandering under $h$, so if $x\in D$, then the $y$-coordinate of $h^n(x)$ goes to $-\infty$ as $n \rightarrow +\infty$ and $H^n(\pi (x))$ goes towards $S$.

For a band $D_{\infty}$ let $\hat D_{\infty}$ be the closure of the component of ${\bf R}^2 \setminus D_{\infty}$ which contains $h(D_{\infty})$.
Let then $E_N = \bigcup \{ \pi ( \hat D_{\infty}) \vert D \ \hbox{of type} \rightarrow \! N \}$ and  $E_S = \bigcup \{ \pi ( \hat D_{\infty}) \vert D \ \hbox{of type} \rightarrow \! S \}$ where $D_{\infty}$ is the band associated to $D$. These sets are closed according to the current hypothesis and are disjoint according to the above remark. Therefore there exist an essential closed curve in $A$ which avoid $E_N$ and $E_S$ and that one can cover by a finite number of elements of $\cal D$ (since we arranged that $\pi ({\cal D})$ covers $A$). These elements are all of type $\leftarrow \! N$ or $\leftarrow \! S$ and therefore, according to the remark, all of the same type. Applying Lemma 4.2 to the Brouwer lines associated to these discs and to all their images under powers of $t$ gives an essential simple closed curve free under $H$ in the annulus.

\medskip

This finishes the proof of the Affirmation and of the Theorem.

\bigskip

{\bf 4.10 Lemma}. {\it Let $(F_n)_{n \geq 1}$ be a numerable sequence of non empty finite sets and $R_i \subset F_i \times F_{i-1}$ a sequence of relations such that for all $f_i \in F_i$ there exist $f_{i-1} \in F_{i-1}$ such that $(f_i, f_{i-1}) \in R_i$. Then there exists a sequence $(f_n)_{n\geq 1}$ such that $f_i \in F_i$ and $(f_i, f_{i-1}) \in R_i$ for all $i\geq 1$.}

\medskip

{\it Proof}. Endow each $F_i$ with the discrete topology. Then $\Pi _{i\geq 1}F_i$ is compact and, if $k\geq 1$, $A_k = \{ (f_1, \ldots , f_k) \vert (f_i, f_{i-1}) \in R_i , 2\leq i \leq k \} \times \Pi _{i>k}F_i$ is a non empty closed subset of  $\Pi _{i\geq 1}F_i$. Therefore $\bigcap _{k>0}A_k \not = \emptyset$. 

\vskip 2cm

\noindent{\bf 5. Winkelnkemper's conjecture}.

\bigskip

In [W], Winkelnkemper proved the following extension of the Poincar\'e-Birkhoff theorem: 

\bigskip

{\bf Theorem}. Let $H$ be a homeomorphism of the closed annulus $S^1 \times [0, 1]$ homotopic to the identity and let $h$ be a lift of $H$ to ${\bf R} \times [0, 1]$ such that $h$ is not conjugate to the translation $t : (x, y) \mapsto (x+1, y)$. Then $h$ has a fixed point or there exists in the annulus an essential simple closed curve free under $H$.

\medskip

Since the hypothesis $h$ not conjugate to $t$ makes sense on the open annulus, he could end his note with the following conjecture. {\it Let $H$ be an orientation and area preserving homeomorphism of the $2$-sphere $S^2$ with two distincts stable fixed points $N$ and $S$, then, if $h$ is a lift of $H$ to the universal cover of $S^2 \setminus \{ N, S \}$ without fixed point, $h$ is conjugate to $t$}. The following theorem is a slightly stronger version of Winkelkemper's conjecture.

\bigskip

{\bf Definition}. A fixed point $p$ of a homeomorphism $H$ is {\bf stable} if there is a basis $\{ U_n \}_{n \in {\bf N}}$ of neighborhoods of $p$ such that $H(U_n) = U_n$ for all $n \in {\bf N}$.

\medskip

{\bf Remarks}. 1) If $H$ has no wandering point, for the fixed point $p$ to be stable, it is enough that there exists a basis $\{ U_n \}_{n \in {\bf N}}$ of neighborhoods of $p$ such that $H(U_n) \subset U_n$. Indeed, if $H(U) \subset U$ for some open set $U$, then $U \setminus H (\overline U)$ is a wandering open set.

2) We can suppose each $U_n$ connected by replacing each $U_n$ by the connected component of $U_n$ which contains $p$.

\bigskip

{\bf Theorem}. {\it Let $H$ be an orientation preserving homeomorphism of $S^2$ without wandering point and with two distinct stable fixed points $N$ and $S$. Then, if $h$ is a lift of $H$ to the universal cover of $S^2 \setminus \{ N, S \}$ without fixed point, $h$ is conjugate to $t$}.

\medskip

{\it Proof}. Let $L$ be a Brouwer line for $h$ which projects down to the interior of a simple arc $\gamma$ in $S^2$ going from $S$ to $N$ as given by Theorem 4.5. Define $W= \bigcap _{n \in {\bf Z}}D(h^n (L))$ and $W'= \bigcap _{n \in {\bf Z}}G(h^n (L))$; these are $h$-invariant closed sets. To show that $h$ is conjugate to $t$ it is enough to show that $W$ and $W'$ are empty. By exchanging $h$ and $h^{-1}$ we can restrict to consider $W$.

Let $\{ U_n \}_{n \in {\bf N}}$ and $\{ V_n \}_{n \in {\bf N}}$ be basis of connected invariant neighborhoods of $N$ and $S$ respectively. Suppose there exists $p \in W$ and choose $n$ such that $p \in \pi ^{-1}(S^2 \setminus U_n \bigcup V_n)$ and a component $\delta$ of $\gamma \bigcap (S^2 \setminus U_n \bigcup V_n)$ joining Fr$U_n$ and Fr$V_n$. Choose then a component $\alpha$ of $\pi ^{-1}(\delta )$, this is a free sub-arc of $L$ joining Fr$\pi ^{-1}(U_n)$ and Fr$\pi ^{-1}(V_n)$ in the $h$-invariant closed set ${\bf R}^2 \setminus \pi ^{-1}(U_n \bigcup V_n)$ and this set is separated by $\alpha $. We can now repeat Winkelnkemper argument: if $k \in {\bf Z}$ is such that $p$ lies between $\alpha$ and $t^k(\alpha )$, then for all $n \geq 0$, $h^{-n}(p)$ lies between $\alpha$ and $h^{-n}t^k(\alpha ) = t^kh^{-n}(\alpha )$ and so, between $\alpha$ and $t^k(\alpha)$. This is a contradiction, since all points of the plane are wandering under the Brouwer homeomorphism $h$.

 \vskip 2cm
{
\centerline {REFERENCES}\bigskip
{ \parindent=5mm 

\item{[A]}   S. Addas-Zanata, {\it Some extensions of the Poincar\'e-Birkhoff theorem to the cylinder and a remark on mappings of the torus homotopic to Dehn twists}, Nonlinearity 18 (2005), 2243-2260.

\medskip

\item{[B]}   L.E.J. Brouwer, {\it Beweis des ebenen Translationssatzes}, Math.
Ann., 72 (1912), 37-54.

\medskip

\item{[BCL]}  F. B\'eguin, S. Crovisier, F. Le Roux, {\it Pseudo-rotations of the open annulus}, preprint arXiv:math.DS/0506041 v1 2Jun 2005.

\medskip

\item{[BH]} M. Bestvina, M. Handel, {\it An area preserving homeomorphism of $T^2$ that is fixed point free but does not move any essential simple closed curve off itself}, Ergod. Th. and Dynam. Sys., 12 (1992), 673-676.

\medskip

\item{[F]}    J. Franks, {\it Generalizations of the Poincar\'e-Birkhoff theorem}, Ann.
of Math., 128 (1988), 139-151.

\medskip

\item{[G]}  L. Guillou, {\it Th\'eor\`eme de translation plane de Brouwer et
g\'en\'eralisations du th\'eor\`eme de Poincar\'e-Birkhoff}, Topology, 33 (1994), 331-351.

\medskip

\item{[K1]}  B.de K\'er\'ekjart\'o, {\it Vorles\" ungen uber Topologie}, Springer, Berlin (1923). 

\medskip

\item{[K2]}  B.de K\'er\'ekjart\'o, {\it The plane translation theorem of Brouwer and the last
geometric theorem of Poincar\'e}, Acta Sci. Math. Szeged, 4 (1928-29), 86-102. 

\medskip

\item{[L]} P. Le Calvez, {\it Une version feuillet\'ee \'equivariante du th\'eor\`eme de translation de Brouwer}, Publications IHES, 102 (2005), 1-98.

\medskip

\item{[S]} A. Sauzet, {\it Application des d\'ecompositions libres \`a l'\'etude des hom\'eomorphismes de surfaces}, PhD thesis, Universit\'e Paris-Nord (2001).

\medskip
  
\item{[T]}  H. Terasaka, {\it Ein Beweis des Brouwerschen ebenen Translationssatzes},
Japan J. of Math., 7 (1930), 61-69. 

\medskip

\item{[W]} H.E. Winkelnkemper, {\it A generalisation of the Poincar\'e-Birkhoff theorem}, Proc. AMS, 102 (1988), 1028-1030.

\par}\bigskip \bigskip \bigskip
{\parindent=30mm
\item{Lucien GUILLOU}  Universit\'e Grenoble 1, Institut Fourier B.P. 74, Saint-Martin-d'H\`eres
38402 (cedex) France 
\par}
\noindent{lguillou@ujf-grenoble.fr}

\end